\documentclass[11pt]{amsart}
\usepackage{url,upref,setspace}
\usepackage{amssymb,amscd,amsthm,amsmath}
\newtheorem{theorem}{Theorem}[section]

\newtheorem{corollary}[theorem]{Corollary}

\newtheorem*{definition}{Definition}

\newtheorem{lemma}[theorem]{Lemma}

\begin{document}
\title{Abelian Groups, Homomorphisms and Central Automorphisms of
  Nilpotent Groups}
\author{Ayan Mahalanobis}
\address{Department of Mathematical Sciences,
  Stevens Institute of Technology, Hoboken, NJ
  07030, USA.}
\thanks{Ayan.Mahalanobis@stevens.edu}
\date{}
\maketitle
\section{Introduction}
It is natural to try to find a necessary and sufficient condition for
a finite group to have an abelian central automorphism group. In this
paper we find a
necessary and sufficient condition in case when the group is finite and
nilpotent, see Theorem \ref{lastthm}. Since a nilpotent group is the
direct product of its Sylow 
subgroups, finding a necessary and sufficient condition for a
nilpotent group is equivalent to finding a necessary and sufficient
condition for a $p$-group. So from now on we work with $p$-groups.
Also, we saw in \cite{ayan2,ayan1} that a $p$-group with a non-trivial abelian 
subgroup of its automorphism 
group can be used to build a key exchange protocol, useful in public
key cryptography. This author in \cite{ayan2,ayan1} used a family of
groups with commutative central automorphism group in a Diffie-Hellman type
key exchange protocol. 

The most notable of the recent publications in the direction of understanding
the central automorphisms of a finite $p$-group are \cite{curran2,curran1,fournelle1,jamali1,curran}. Jamali and Mousavi in 
\cite{jamali1} provide a necessary and sufficient condition for a
$p$-group $G$ of class 2, for an odd prime $p$, to have an elementary abelian
central automorphism group and Curran \cite{curran} studies groups
which have an
abelian automorphism group, i.e., Miller groups. In a Miller group all
automorphisms are central, so our work contributes
to the study of the central automorphisms of nilpotent groups as
well as the study of Miller groups. For an introduction to Miller
groups see \cite{earnley,morigi} or \cite[Section 2.7]{ayan1}.
\begin{definition}[PN Group]
A group $G$ is a purely non-abelian group if it doesn't have any nontrivial
abelian direct factor.     
\end{definition}
Adney and Yen in
\cite[Theorem 4]{adney} proved a necessary and sufficient condition for a
finite PN $p$-group of class 2 to have an abelian central automorphism
group. In this 
article we extend that result and to a certain extent their argument
to arbitrary 
finite PN $p$-groups. We split the problem into three parts:
\begin{itemize}
\item[a] Reduce the problem to a problem about abelian groups and
  homomorphisms between abelian groups.
\item[b] Solve the problem in finite abelian groups.
\item[c] Bring the solution from abelian groups back to non-abelian
  finite PN $p$-groups.
\end{itemize}
This is not the first time that the theory of abelian groups has been used in
understanding the automorphisms of non-abelian $p$-groups. Sanders
\cite{sanders} used a technique quite similar to ours to count the
number of central automorphisms in a finite PN (purely non-abelian)
$p$-group. 
\section{Central Automorphisms}
Let $\text{Aut}_c(G)$ be the group of central automorphisms of a group $G$. An
automorphism $\sigma\in\text{Aut}(G)$ is called a central automorphism if
$g^{-1}\sigma(g)\in Z(G)$ for all $g\in G$, or equivalently, central
automorphisms are the centralizer of the group of inner automorphisms.

There is another way to think about the central automorphisms. Let
$\sigma\in\text{Aut}_c(G)$, then corresponding to $g\in G$ there is a
$z_{g,\sigma}\in Z(G)$ such that $\sigma(g)=gz_{g,\sigma}$. Corresponding to
$\sigma\in\text{Aut}_c(G)$ one can define a map $\phi_\sigma:G\rightarrow
Z(G)$ as follows: \[\phi_\sigma(g)=z_{g,\sigma}\]

It is straightforward to show that the map $\phi_\sigma$ is a
homomorphism. Hence corresponding to $\sigma\in\text{Aut}_c(G)$ there is
$\phi_{\sigma}\in\text{Hom}(G,Z(G))$. It is known that for  PN groups the
converse is true, see \cite[Theorem 1]{adney}. 

There is a connection \cite[Theorem 3]{adney} between commutativity of
the group 
of central automorphisms $\text{Aut}_c(G)$ and commutativity of the
homomorphisms $\phi_\sigma$. 

Assume that $\text{Aut}_c(G)$ is commutative, then for two maps
$\tau, \sigma\in\text{Aut}_c(G)$ we have that
$\tau(\sigma(g))=\sigma(\tau(g))$ which is the same as
$$g\phi_\sigma(g)\phi_\tau(g)\phi_\sigma(\phi_\tau(g))=g\phi_\tau(g)\phi_\sigma(g)\phi_\tau(\phi_\sigma(g))$$ 
implying, 
$\sigma,\tau\in\text{Aut}_c(G)$ commute if and only if
$\phi_\sigma,\phi_\tau\in\text{Hom}(G,Z(G))$ commute. 
Notice that since $Z(G)$ is an abelian group hence
$\phi_{\sigma}(G^\prime)=1$. So 
corresponding to $\phi_{\sigma}:G\rightarrow Z(G)$ one can define
$\phi_{\sigma}^\prime:\dfrac{G}{G^\prime}\rightarrow Z(G)$ as
$\phi_{\sigma}^\prime(xG^\prime)=\phi_{\sigma}(x)$. Clearly
$\phi_{\sigma}^\prime$ is a homomorphism. 

Consider the map $\lambda:Z(G)\rightarrow G/G^\prime$ given by the diagram
\begin{equation}\label{lambda}
\begin{CD}
Z(G) @>\iota>>G @>\pi>>G/G^\prime
\end{CD}
\end{equation}
where $\iota$ and $\pi$ are the inclusion and the natural surjection
respectively. 
\begin{theorem}
Let $\sigma,\tau\in\text{Aut}_c(G)$ then $\sigma\circ\tau=\tau\circ\sigma$ if
and only if
$\phi_\sigma^\prime\circ\lambda\circ\phi_\tau^\prime=\phi_\tau^\prime\circ\lambda\circ\phi_\sigma^\prime$.  
\end{theorem}
\begin{proof}
The proof follows from the above discussion.
\end{proof}
This theorem enables us to think about commutativity of the group of central
automorphisms of a nonabelian group $G$ in terms of abelian groups
$\dfrac{G}{G^\prime}$ and $Z(G)$ and homomorphisms between these
abelian groups. So the 
problem ``when is $\text{Aut}_c(G)$ commutative?'' for a PN group is
transformed  into a problem involving abelian groups and homomorphisms between
abelian groups. 

This also enables us to ask more general questions about abelian groups and
homomorphisms between abelian groups
that is the object of our study in the next section.
\section{c-maps}
\begin{definition}[c-map]
Let $A$ and $B$ be abelian groups and $\lambda:A\rightarrow B$ a
homomorphism. Then $\lambda$ is a c-map if $f\lambda g=g\lambda f$ for all
$f,g\in\text{Hom}(B,A)$. The set of all c-maps form a subgroup of
$\text{Hom}(A,B)$. We call $\lambda$ a trivial c-map if $f\lambda g=0$ for all
$f,g\in\text{Hom}(B,A)$.  
\end{definition} 
In this section we investigate necessary and sufficient conditions for
two finite abelian groups $A$ and $B$ and $\lambda:A\rightarrow B$ a
homomorphism between them to be a 
c-map. Of course if $\lambda\equiv 0$ then trivially $f\lambda g=0$ for any
$f,g\in\text{Hom}(B,A)$ and $\lambda$ is a c-map. 

For the rest of this section we fix $A$ and $B$ to be two finite nonzero
abelian $p$-groups. From the fundamental theorem of finite abelian groups, we
have 
\[A=A_1\oplus A_2\oplus\cdots\oplus A_{n}\;\;\;\; n\geq 1\]
\[B=B_1\oplus B_2\oplus\cdots\oplus B_{m}\;\;\;\; m\geq 1\]
where 
$$A_i=\langle
a_i\rangle,\;\;\text{exp}(A_i)=p^{\alpha_i}\;\;\;\text{and}\;\;\;
a=\alpha_1\geq\alpha_2\geq\ldots\geq\alpha_n$$
$$B_i=\langle b_i\rangle,\;\;\text{exp}(B_i)=p^{\beta_i}\;\;\;\text{and}\;\;\;b=\beta_1\geq\beta_2\geq\ldots\geq\beta_m$$ are
decompositions of $A$ and $B$ as direct sums of cyclic $p$-groups.
For the rest of this article we fix $A_1,A_2,\ldots,A_{n}$ as a
fixed decomposition of $A$ and $B_1,B_2,\ldots,B_{m}$ as a fixed
decomposition of $B$ and $a_i$ a fixed generator for $A_i$, $1\leq
i\leq n$ and $b_i$ a fixed generator for $B_i$, $1\leq i\leq m$.

If $\lambda:A\rightarrow B$ is a homomorphism, we can write
$\lambda=\lambda_1+\lambda_2+\cdots+\lambda_m$ where $\lambda_i:A\rightarrow
B_i$ is defined as follows:
\[
\lambda_i(a)=\pi_i\lambda(a)
\]
where $\pi_i:B\rightarrow B_i$ is the projection.
We can abuse the notation a little bit and consider
$\lambda_i:A\rightarrow B$. Of course, one can formalize
this trivially. It is easy to see that if each $\lambda_i$ is a c-map, then
$\lambda$ is a c-map.  

Let
$\mathcal{R}=\{x\in A\;|\; |x|\leq p^b\}$. 
It is known that 
\[\mathcal{R}=\sum\limits_{f\in\text{Hom}(B,A)}f(B)\]
and since $\mathcal{R}$ is an abelian group contained in $A$ and if
$x_1+x_2+\ldots +x_n\in\mathcal{R}$ then from the definition of
$\mathcal{R}$ it follows that $x_i\in\mathcal{R}$ for $i=1,2,\ldots,n$. Hence 
$\mathcal{R}=\mathcal{R}_1\oplus\mathcal{R}_2\oplus\ldots\oplus\mathcal{R}_n$
where $\mathcal{R}_i=\langle r_i\rangle$ and $\mathcal{R}_i=\mathcal{R}\cap
  A_i$ for each $i$. We assume that
$\text{exp}(\mathcal{R}_1)=p^{n_1}$ and
$\text{exp}(\mathcal{R}_2\oplus\mathcal{R}_3\oplus\ldots\oplus\mathcal{R}_n)=p^{n_2}$.
Clearly $a\geq n_1\geq n_2$.

Let $e_{ij}:B\rightarrow A$ be defined as
\begin{equation}\label{equation1}
e_{ij}(b_k)=\left\{
\begin{array}{lll}
p^{\max(0,\alpha_i-\beta_j)}a_i &\text{if}\;\; j=k\\
0 &\text{otherwise}
\end{array}\right.
\end{equation}
It follows from \cite[Section 5.8]{scott} that
$\{e_{ij}\}$, $i=1,2,\ldots,n$; $j=1,2,\ldots,m$ is a basis for
$\text{Hom}(B,A)$ under addition. 

From Equation \ref{equation1} either $r_i=a_i$ or 
$r_i=p^{\alpha_i-\beta_j}a_i$. Since $\beta_1$ is the maximum possible
hence $\alpha_i-\beta_1$ is the least possible for a fixed $i$ and for
all $j\;(1\leq j\leq m)$. From this we conclude that
$e_{i1}(b_1)=r_i$ for all $i$.

We state some easy and well known facts in the following lemma whose
proof follows from the above discussion.
\begin{lemma}[]
\begin{itemize}
\item[(i)] $\mathcal{R}=A[p^b]$.
\item[(ii)]
  $\mathcal{R}=\mathcal{R}_1\oplus\mathcal{R}_2\oplus\ldots\oplus\mathcal{R}_n$ where $\mathcal{R}_i=\mathcal{R}\cap A_i$.
\item[(iii)] If $\mathcal{R}_i=\langle r_i\rangle$, then $r_i=a_i$ if
  $b\geq \alpha_i$ or $r_i=p^{\alpha_i-\beta_1}a_i$ if $b<\alpha_i$.
\item[(iv)] $e_{i1}(b_1)=r_i$ for all $i$.
\item[(v)] $e_{1j}(b_j)=p^{\text{max}(0,n_1-\beta_j)}r_1$ for $j>1$.
\end{itemize}
\end{lemma} 
\begin{theorem}
A homomorphism $\lambda:A\rightarrow B$ is a c-map if and only if
$e_{ij}\lambda e_{kl}=0$ whenever $i\neq k$ or $j\neq l$.
\end{theorem}
\begin{proof}
Let $f\lambda g=g\lambda f$ for all $f,g\in\text{Hom}(B,A)$ then $e_{ij}\lambda
e_{kl}=e_{kl}\lambda e_{ij}$. If $i\neq k$ then $e_{ij}\lambda
e_{kl}=0$. Since the image of $e_{ij}\lambda e_{kl}$ is in $A_i\cap A_k$ where
$A_i\cap A_k=0$ whenever $i\neq k$. If $j\neq l$ then $e_{kl}\lambda
e_{ij}(b_l)=0$. Hence, $e_{ij}\lambda e_{kl}=0$ whenever $i\neq k$ or
$j\neq l$. 

Conversely, notice that any $f\in\text{Hom}(B,A)$ can be written as
$f=\sum\limits_{j=1}^m\sum\limits_{i=1}^nn_{ij}e_{ij}$ where
$n_{ij}\in\mathbb{Z}$. Hence we can write
$f=\sum\limits_{j=1}^m\sum\limits_{i=1}^nn_{ij}e_{ij}$ where
$n_{ij}\in\mathbb{Z}$ and
$g=\sum\limits_{j=1}^m\sum\limits_{i=1}^nn_{ij}^\prime e_{ij}$ where
$n_{ij}^\prime\in\mathbb{Z}$.
Hence 
\begin{eqnarray*}
\lefteqn{f\lambda
  g}\\
&=&\left(\sum\limits_{j=1}^m\sum\limits_{i=1}^nn_{ij}e_{ij}\right)\lambda\left(\sum\limits_{k=1}^m\sum\limits_{l=1}^nn_{kl}^\prime
  e_{kl}\right)\\
&=&\sum\limits_{j=1}^m\sum\limits_{i=1}^n\left(n_{ij}e_{ij}\lambda
  n_{ij}^\prime e_{ij}\right)\;\;\;\text{since}\;\; e_{ij}\lambda
  e_{kl}=0\;\;\text{whenever}\;\;i\neq k\;\text{or}\;j\neq l\\
&=&\sum\limits_{j=1}^m\sum\limits_{i=1}^nn_{ij}n_{ij}^\prime\left(e_{ij}\lambda
  e_{ij}\right)\\
&=&\sum\limits_{j=1}^m\sum\limits_{i=1}^nn_{ij}^\prime
  n_{ij}\left(e_{ij}\lambda e_{ij}\right)=g\lambda f.
\end{eqnarray*}
\end{proof}
We use the above theorem to prove: 
\begin{theorem}\label{thm1}
A homomorphism $\lambda:A\rightarrow B$ is a c-map if and only if 
\begin{eqnarray}
\label{one1}\lambda(r_u)&\in&p^{n_1}B\;\;\;\;\;\;\;\;u>1\\
\label{one2}\lambda_1(r_1)&\in&\langle
p^{k^\prime}b_1\rangle\;\;\;\;\;\;k^\prime=\min(n_1,\max(n_2,\beta_2))\\ 
\label{one3}\lambda_j(r_1)&\in&\langle p^{n_1}b_j\rangle\;\;\;\;\;j>1.
\end{eqnarray}
Moreover, if $\lambda$ satisfies the above conditions and
$\langle\lambda_1(r_1)\rangle =\langle p^kb_1\rangle$ where $k^\prime\leq k<n_1$, then $\lambda(\mathcal{R})=p^kB_1$. 
\end{theorem}
\begin{proof}
We assume that conditions (\ref{one1}), (\ref{one2}) and (\ref{one3}) are
satisfied. If $u>1$ then $\lambda(r_u)\in p^{n_1}B$.
Hence $e_{st}\lambda e_{uv}=0$ for all $s,t$ and $v$. If $u=1$ then for $s>1$,
$e_{st}\lambda e_{uv}(b_v)\in \mathcal{R}_s$. Now since
$\text{exp}(\mathcal{R}_s)\leq p^{n_2}$ and 
$\lambda_j(r_1)\in\langle p^{n_2}b_j\rangle$ for $1\leq j\leq m$,
hence we have $e_{st}\lambda e_{uv}=0$. From the earlier discussion it
follows that 
$e_{11}(b_1)=r_1$. Now notice that for $t>1$, 
$e_{1t}\lambda e_{11}(b_1)=e_{1t}\lambda(r_1)=0$ and $e_{11}\lambda
e_{1t}(b_t)=e_{11}\lambda(p^{\max(0,n_1-\beta_t)}r_1)=0$, from the
definition of $k^\prime$.

Conversely, we assume that $e_{st}\lambda e_{uv}=0$ whenever $s\neq u$ or
$t\neq v$. Now $e_{11}(b_1)=r_1$ then for
a $j>1$ $e_{1j}\lambda e_{11}=0$. That says that
\[e_{1j}\lambda_j(r_1)=0\;\;\text{for each}\;\; j>1\] 
Now from Equation \ref{equation1} either $e_{1j}(b_j)=r_1$ or
$e_{1j}(b_j)=p^{n_1-\beta_j}r_1$. In the first case clearly
$\lambda_j(r_1)\in\langle p^{n_1}b_j\rangle$. In the second case it follows
from 
$e_{1j}\lambda(r_1)=0$ that $\lambda_j(r_1)=p^{\beta_j} b_j=0$. In either
case $\lambda_j(r_1)\in\langle p^{n_1}b_j\rangle$. This proves (\ref{one3}).

Pick a $u>1$, then $e_{1j}\lambda e_{u1}=0$ for all $j$. This implies
that $e_{1j}\lambda(r_u)=0$ for all $j$. Since $e_{1j}(b_j)=r_1\;\;
\text{or}\;\; p^{n_1-\beta_j}r_1$ hence $\lambda_j(r_u)\in\langle
p^{n_1}b_j\rangle\;\;\text{or}\\ \lambda_j(r_u)=0$ for $u>1$ and all $j$. This
proves (\ref{one1}).

Choose $u>1$ such that
$\text{exp}(\mathcal{R}_2\oplus\mathcal{R}_3\oplus\cdots\oplus\mathcal{R}_n)=\text{exp}(\mathcal{R}_u)=p^{n_2}$
and then  we have that $e_{u1}(b_1)=r_u$, then clearly
$e_{u1}\lambda e_{11}=0$ implies that $e_{u1}\lambda(r_1)=0$, hence
$\lambda_1(r_1)\in\langle p^{n_2}b_1\rangle$. 

On the other hand $e_{11}\lambda e_{1j}=0$ for $j>1$ implies that

\noindent $e_{11}\lambda(p^{\max(0,n_1-\beta_j)}r_1)=0$ which implies that
$p^{\max(0,n_1-\beta_j)}e_{11}\lambda(r_1)=0$ for all $j>1$. This
gives us that $\lambda_1(r_1)\in\langle
p^{k^\prime}b_1\rangle$. The above two arguments proves (\ref{one2}).  

The later assertion of the theorem clearly follows from the fact that
$\beta_2\leq k$ and hence $p^{n_1}b_j=0$ for all $j>1$.
\end{proof}
\begin{corollary}
If $n_1=n_2$ and $\lambda:A\rightarrow B$ is a c-map then $\lambda$ is
a trivial c-map. 
\end{corollary}
\begin{proof}
It follows from Theorem \ref{thm1}, since $k^\prime=n_1$ that $\lambda(r_i)\in p^{n_1}B$ for all
$i$, hence $f\lambda(r_i)=0$ for all $i$ and hence $f\lambda g=0$.
\end{proof}
\begin{corollary}
If $b=n_1$ and $\lambda$ is a c-map, then
$\lambda(\mathcal{R}_2\oplus\mathcal{R}_3\oplus\cdots\oplus\mathcal{R}_n)=0$.
This automatically implies that $\lambda(\mathcal{R})$ is cyclic.
\end{corollary}
\begin{proof}
From Theorem \ref{thm1} it follows that $\lambda(r_i)=0$ for $i>1$. 
\end{proof}
Now assume there is a $j$ such that 
$\text{exp}(B_j)=\text{exp}(B_1)$ and $j>1$. Then clearly
$e_{1j}(b_j)=r_1$. If $\lambda$ is a c-map then $e_{11}\lambda e_{1j}(b_j)=0$
implying that $e_{11}\lambda(r_1)=0$. This tells us that
$\lambda_1(r_1)=p^{n_1}b_1$ which makes $\lambda$ a trivial c-map. We just
established that a necessary condition for $\lambda$ to be a
nontrivial c-map is that $\text{exp}(B_1)>\text{exp}(B_j)$ for $j>1$.

If $\lambda:A\rightarrow B$ is any homomorphism, let $p^c=\text{exp}(\text{ker}(\lambda))$, clearly $c\leq a$. We now
find a necessary and sufficient condition for
$\text{ker}(\lambda)\subseteq\mathcal{R}$. Recall that $\mathcal{R}=A[p^b]$. 
\begin{lemma}
$\text{ker}(\lambda)\subseteq\mathcal{R}$ if and only if $c\leq n_1$.
\end{lemma}
\begin{proof}
Let $\text{ker}(\lambda)\subseteq\mathcal{R}$. Since
$p^{n_1}\mathcal{R}=0$ hence 
$p^{n_1}\text{ker}(\lambda)=0$. This proves that $c\leq n_1$.
Conversely, assume that $c\leq n_1$. Let $x\in\text{ker}(\lambda)$
then $p^{n_1}x=0$ hence $x\in\mathcal{R}$. 
\end{proof}
The next lemma comes in handy to settle the question: if $c\geq n_1$
then are there any nontrivial c-maps? 
\begin{lemma}
$p^cB=\bigcap\limits_{f:B\rightarrow\text{ker}(\lambda)}\text{ker}(f)$.
\end{lemma}
\begin{proof}
Let $x\in p^cB$ then $x=p^cy$ for some $y\in B$. Then for any
$f\in\text{Hom}(B,\text{ker}(\lambda))$, $f(x)=p^cf(y)=0$, since the exponent
of $\text{ker}(\lambda)$ is $c$. So
$p^cB\subseteq\bigcap\limits_{f:B\rightarrow\text{ker}(\lambda)}\text{ker}(f)$.  
Conversely, assume that $x\not\in p^cB$ then the image of $x$ in
$\dfrac{B}{p^cB}$ is nontrivial and
exp$\left(\dfrac{B}{p^cB}\right)=p^c=\text{exp}(\text{ker}(\lambda))$. Hence
there is a $f\in\text{Hom}(B,\text{ker}(\lambda))$ such that $f(x)\neq 0$. 
\end{proof}
Taking the two previous lemmas together we show that
\begin{lemma}\label{lm1}
If $\lambda:A\rightarrow B$ is a c-map then
$\lambda(\mathcal{R})\subseteq p^cB$. It
follows that if $\lambda$ is a c-map and $c\geq n_1$ then $\lambda$ is a
trivial c-map. 
\end{lemma}
\begin{proof}
Since $\lambda$ is a c-map, hence $F\lambda g=g\lambda F$ where
$F\in\text{Hom}(B,\text{ker}(\lambda))$ and $g\in\text{Hom}(B,A)$. Then
clearly $F\lambda g=0$ for all $g\in\text{Hom}(B,A)$. Hence
$F\lambda(\mathcal{R})=0$ for all
$F\in\text{Hom}(B,\text{ker}(\lambda))$. Hence
$$\lambda(\mathcal{R})\subseteq\bigcap\limits_{F:B\rightarrow\text{ker}(\lambda)}  
\text{ker}(F)=p^cB.$$ The rest of the argument follows from the fact
that $\text{exp}(\mathcal{R})=p^{n_1}$.   
\end{proof}
We just saw that for all intended purposes of understanding c-maps
$c\geq n_1$ is
irrelevant, because then $\lambda$ is a trivial c-map. So, from now on
we will work with $c<n_1$ which implies that
$\text{ker}(\lambda)\subseteq\mathcal{R}$.  

This has little relevance to the flow of arguments towards the proof
of the main result but is of independent interest. Using the same
method as above one can easily prove that  
\[p^aB=\bigcap\limits_{f:B\rightarrow A}\text{ker}f.\]
This yields a lemma, whose proof we leave to the reader and is
corroborated by Theorem \ref{thm1}:
\begin{lemma}
A homomorphism $\lambda$ is a trivial c-map if and only if
$\lambda(\mathcal{R})\subseteq p^{n_1}B$.
\end{lemma}

We are now in a position to prove the main theorem of this article. 
\begin{theorem}\label{thm2}
Let $\lambda:A\rightarrow B$ be a homomorphism. Then $\lambda$ is
a non-trivial c-map if and only if 
\begin{eqnarray*}
&\lambda(\mathcal{R})= p^kB & \text{where}\;\; c\leq k<n_1\\
&\text{and}\;\;\;\dfrac{\mathcal{R}}{\text{ker}(\lambda)}& \text{is cyclic}.
\end{eqnarray*}
\end{theorem}
\begin{proof}
The only if part follows from Theorem \ref{thm1} and Lemma 3.8.

To see the if part, assume that $\lambda(\mathcal{R})=p^kB$ where
$c\leq k<n_1$ and 
$\dfrac{\mathcal{R}}{\text{ker}(\lambda)}$ is cyclic.
Without loss of generality we assume that
$\dfrac{\mathcal{R}}{\text{ker}(\lambda)}\cong\langle
p^kb_1\rangle$. Hence there is some $r\in\mathcal{R}$ such that 
$\lambda(r)=p^kb_1$. Also from $n_1>c$, it follows that
$|r|=n_1$. We show that $f\lambda g(b_i)=0$ for all
$f,g\in\text{Hom}(B,A)$ and $i\geq 2$.

It is clear that $|b_i|\leq p^k$ for $i\geq 2$ and $g(b_i)=sr+u$
where $u\in\text{Ker}\lambda$ and $p^{n_1-k}|s$. Then $s=s^\prime
p^{n_1-k}$. Hence $f\left(\lambda(s^\prime
  p^{n_1-k}r+u)\right)=f\left(s^\prime p^{n_1}b_1\right)=0$ and this proves the theorem.
\end{proof}
It is interesting to note what happens in case of a c-map $\lambda$ such that
$\text{exp}(\text{ker}(\lambda))=\text{exp}(\text{coker}(\lambda))$
which implies that $c\leq b$.

Notice that $p^{b-c}p^cB=p^bB=0$. Now if $\lambda(x)\in p^cB$ then
$\lambda\left(p^{b-c}x\right)=0$, i.e., $p^{b-c}x\in\text{ker}(\lambda)$. This
says that $p^cp^{b-c}x=p^bx=0$ which implies $x\in\mathcal{R}$. Now assume that
$\text{exp}(\text{ker}(\lambda))=\text{exp}(\text{coker}(\lambda))$ then we
have that $p^cB\subseteq\text{Image}(\lambda)$. Hence for any $y\in p^cB$ there
is a $x\in A$ such that $\lambda(x)=y$. This implies that
$x\in\mathcal{R}$. Hence $p^cB\subseteq\lambda(\mathcal{R})$. Using
Lemma \ref{lm1} we just proved the following lemma:
\begin{lemma}
Let $\lambda:A\rightarrow B$ be a c-map and
$\text{exp}(\text{ker}(\lambda))=\text{exp}(\text{coker}(\lambda))$
then $\lambda(\mathcal{R})=p^cB$. 
\end{lemma} 
\section{Back to $p$-groups}
In this section we use the theorems from the last section to find a necessary
and sufficient condition, in the same spirit as in Adney and Yen
\cite[Theorem 4]{adney}, for the group of central automorphisms of a finite
PN p-group $G$ to be 
abelian. We have seen before that there is a one-to-one correspondence between
the central automorphisms in $G$ and homomorphisms from
$\dfrac{G}{G^\prime}$ to $Z(G)$. Now the central automorphisms commute if and
only if $\lambda: Z(G)\rightarrow\dfrac{G}{G^\prime}$ defined by
$\lambda(x)=xG^\prime$ is a c-map. We use all the notation from the last
section with the understanding that $Z(G)$ represents $A$ and
$\dfrac{G}{G^\prime}$ represents $B$. Since the group $G$ is no longer
abelian, even though $\dfrac{G}{G^\prime}$ and $Z(G)$ are abelian, we will no
longer use $+$ to denote the group operation. Clearly then the kernel of
$\lambda$ is 
$Z(G)\cap G^\prime$. Notice that for a $p$-group $G$ of class 2,
$G^\prime\subseteq Z(G)$ and
$\text{exp}(G^\prime)=\text{exp}\left(\dfrac{G}{Z(G)}\right)$. This
means that
$\text{exp}(G^\prime)\leq\text{exp}\left(\dfrac{G}{G^\prime}\right)$.
This clearly implies that $c\leq b$.    
\begin{theorem}\label{lastthm}
Let $G$ be a PN p-group and $p^c=\text{exp}(Z(G)\cap G^\prime)$. Then
the central 
automorphisms of $G$ commute if and only if either
\begin{eqnarray*}
&\lambda(\mathcal{R})\subseteq
\left(\dfrac{G}{G^\prime}\right)^{p^{n_1}}&\\
\text{or}\\
&\lambda(\mathcal{R})= \left(\dfrac{G}{G^\prime}\right)^{p^k} &
\text{where}\;\; c\leq k<n_1\\ 
&\text{and}\;\;\;\dfrac{\mathcal{R}}{Z(G)\cap G^\prime}&\text{is cyclic}. 
\end{eqnarray*}
\end{theorem}
\begin{proof}
This theorem follows from Theorem \ref{thm2}. Notice that $\lambda$ in this
case is the map from Equation \ref{lambda}.
\end{proof}

We should mention the relation of our theorem with that of the
\cite[Theorem 4]{adney}. There the authors work only with $p$-groups of
class 2. In that case we have that $\text{ker}(\lambda)=\text{coker}(\lambda)$
and hence $\lambda(\mathcal{R})=p^cB$ which is the same as Adney and Yen's
condition $\mathcal{R}=\mathcal{K}$. Again since in a $p$-group of
class 2, $\text{ker}(\lambda)=G^\prime$
their condition reads like $\dfrac{\mathcal{R}}{G^\prime}=\langle
x_1^{p^c}G^\prime\rangle$.

\noindent\textbf{Acknowledgement:}
The author wishes to thank Fred Richman for his help and guidance in
preparation of this manuscript. The author is indebted to the referee
for his comments and suggestions which has helped in a better
presentation of this article.
\nocite{fuchs,kaplansky,GAP4,earnley,curran,jamali1,morigi,curran1,fournelle1,curran2,curran3,fournelle}
\bibliography{signature_scheme}
\bibliographystyle{abbrv}
\end{document}